\documentclass[11 pt]{amsart}
\usepackage{graphicx} 
\usepackage{amsmath}
\usepackage{mathtools}
\usepackage
{amstext, amscd, setspace, tikz-cd,mathtools, enumerate, graphics, latexsym, siunitx}
\usepackage{pst-node}
\usepackage{tikz-cd} 
\usepackage{amssymb}
\usepackage{hyperref}

\usepackage{PTSerif} 
\usepackage{graphicx}
\usepackage[cmtip,all]{xy}

\newcommand{\R}{{\mathbb R}}

\newcommand{\C}{{\mathbb C}}

\newcommand{\Hom}{{\rm Hom}}

\newcommand{\Tr}{{\rm Trace~}}

\input{xypic}
\xyoption{all}
\newtheorem{theorem}{Theorem}[section]

\newtheorem{defn}[subsubsection]{{\bf Definition}}
\newtheorem{lemma}[theorem]{Lemma}

\newtheorem{thm}[subsubsection]{{\bf Theorem}}

\newtheorem{rem}[subsubsection]{{\bf Remark}}
\theoremstyle{definition}

\newtheorem{remark}{Remark}[section]
\numberwithin{equation}{section}

\title[Representation equivalence of Lattices]{Representation equivalence of Lattices in Lie Groups}

\author{Chandrasheel Bhagwat and Kaustabh Mondal}

\address{Indian Institute of Science Education and Research, Dr.\,Homi Bhabha Road, Pashan, Pune 411008,  INDIA.}
\email{cbhagwat@iiserpune.ac.in, \ kaustabh.mondal@students.iiserpune.ac.in}

\date{\today}
\subjclass[2020]{22E40, 22E45, 53C35}

\keywords{Lattices in Lie Groups, Representation Equivalence, Trace Formula}

\begin{document}
\begin{abstract} Let $\Gamma_1$ and $\Gamma_2$ be two lattices of finite covolume in a semisimple Lie group $G$. We prove a spectral rigidity result for the representation spectra of the right regular representations $L^2(\Gamma_1 \backslash G)$ and $L^2(\Gamma_2 \backslash G)$ of $G$. This can be thought of as an analogue of the strong multiplicity one theorem 
and it generalises \cite[Thm 1.1]{BR} by the first author and Rajan to the case of non-uniform lattices.



\end{abstract}

\maketitle
\section{Introduction} 

\subsection{Representation Equivalence of Lattices}

The relation between the arithmetic properties of discrete subgroups of Lie groups and the spectrum of associated  locally symmetric spaces is well known. 



\smallskip

Let $G$ be a semisimple Lie group of real rank one. Let $\Gamma$ be a lattice in $G$ such that the covolume $\text{vol}\ (\Gamma \backslash G)$ (with respect to the unique $G$-invariant measure on $\Gamma \backslash G$) is finite. Consider the Hilbert space 
$$ L^2(\Gamma \backslash G)=\{ \phi : G \longrightarrow \C : \phi(\gamma x) =\phi(x),\ \int\limits_G |\phi(x)|^2\ dx < \infty\}.$$
\noindent $G$ acts on $L^2(\Gamma \backslash G)$ via $g.\phi (x)=\phi(xg)$ for all $x, g \in G$. This gives rise to a right regular representation $R_{\Gamma}$ of $G$ on $L^2(\Gamma \backslash G)$. Let $L^2_d(\Gamma \backslash G)$ be the maximal semisimple subrepresentation of $L^2(\Gamma \backslash G)$ (\textit{discrete spectrum}). Then $$L^2_d(\Gamma \backslash G)\cong \widehat{\bigoplus_{\pi \in \widehat{G}}}m(\pi, \Gamma)\ \pi,$$
where, the finite multiplicities $m(\pi, \Gamma)$ can be interpreted as $\text{dim}_{\C}\ \Hom_G(\pi, L^2(\Gamma \backslash G))$.

\noindent The orthogonal complement of $L^2_d(\Gamma \backslash G)$ in $L^2(\Gamma \backslash G)$ is called \textit{continuous spectrum} and is denoted by $L^2_c(\Gamma \backslash G)$. It admits a Hilbert direct integral decomposition with respect to a certain Radon measure on the unitary dual of $G$ (see Sec. \ref{sec: - Spectral Decomposition}\ ).

\begin{defn}
Two lattices $\Gamma_1$ and $\Gamma_2$ are called representation equivalent if the right regular representations $(R_{\Gamma_1}, L^2(\Gamma_1 \backslash G))$ and $(R_{\Gamma_2}, L^2(\Gamma_2 \backslash G))$ are isomorphic as representations of $G$. 
\end{defn}

\noindent The purpose of this short note is to provide a sufficient condition for the representation equivalence of two lattices $\Gamma_1$ and $\Gamma_2$. This can be thought of as an analogue of the classical strong multiplicity one theorem for the group ${\rm GL}(n)$. Before we describe our main result, let us define some notations.\medskip

Since $G$ is of real rank one, every proper parabolic subgroup is minimal. A parabolic subgroup $P$ of $G$ has a Levi decomposition $LN$, where, $L$ is a Levi part of $P$ and $N$ is the unipotent radical of $P$. Moreover, $P$ possesses a decomposition (called Langlands decomposition) as $P=MAN$, where, $A$ is the split torus in $L$ and $M$ is a compact subgroup contained in a maximal compact subgroup $K$; in fact, $M=Z_K(A)$ (see Subsec. \ref{Associate parabolics in $G$} for details).

\begin{defn}
    For a lattice $\Gamma$ in $G$, a proper parabolic subgroup $P$ is called $\Gamma$-cuspidal if $P \cap \Gamma \subset MN$; $\dfrac{N}{N \cap \Gamma}$ and $\dfrac{MN}{MN \cap \Gamma}$ are compact. 
    \end{defn}\smallskip

\noindent The set of all $\Gamma$-cuspidal parabolic subgroups is denoted by $E(G, \Gamma)$. Then $\Gamma$ acts on $E(G , \Gamma)$ by conjugation, and the set of cusps for $\Gamma$ is parametrised by the orbits of this action. Note that if $\Gamma$ is a cocompact lattice then $E(G , \Gamma)$ is empty i.e., there is no cusp for $\Gamma$. \smallskip

\begin{defn}
    Two lattices $\Gamma_1$ and $\Gamma_2$ of finite covolume in $G$ are said to be isocuspidal if there exists a finite set $\{ P_1, P_2, \dots, P_n\}$ of cuspidal parabolic subgroups which serves as a common set of representatives of cusps for both lattices $\Gamma_1$ and $\Gamma_2$.
\end{defn}

\smallskip
\begin{rem}
    If $\Gamma_1$ and $\Gamma_2$ both have only one cusp then they are isocuspidal. 
\end{rem}
\smallskip

The following {\it Spectral Rigidity Theorem} is the main result of this article.\smallskip

\begin{thm} \label{main thm}  

Let $\Gamma_1$ and $\Gamma_2$ be two lattices with finite covolume in $G$ such that the following holds:\medskip

\begin{enumerate}
\item The lattices $\Gamma_1$ and $\Gamma_2$ are isocuspidal and there exists a parabolic subgroup $P = MAN$ in the common set of  representatives of cusps for $\Gamma_1$ and $\Gamma_2$ such that
$L^2((\Gamma_{1} \cap M) \backslash M) \cong L^2( (\Gamma_{2} \cap M) \backslash M)$ as representations of $M$.\medskip

\item The multiplicities $m(\pi, \Gamma_1) = m(\pi, \Gamma_2)$ for all but finitely many $\pi \in \widehat G$.

\end{enumerate}\medskip

\noindent Then the lattices $\Gamma_1$ and $\Gamma_2$ are representation equivalent in $G$.

\end{thm}

\begin{rem}
Note that if $\Gamma_1$ and $\Gamma_2$ are cocompact lattices then the condition (1) in Thm. \ref{main thm} is vacuously true and hence it  essentially recovers \cite[Thm 1.1]{BR}. However, in the case of non-cocompact lattices, condition (1) assures the equality of the continuous spectra (see Subsec. \ref{subsec: Equality of continuous Spectrum} for details).
\end{rem}

\begin{rem}
    If $G$ is reductive but not semisimple, then for congruence lattices $\Gamma_1$ and $\Gamma_2$ in $G$, it is expected that only condition (2) as in Thm. \ref{main thm} is sufficient to conclude $L^2_d(\Gamma_1 \backslash G) \cong L^2_d(\Gamma_2 \backslash G)$. The proof of this will appear in Ph.D thesis of second author (\cite{Km}). 
\end{rem}

\subsection{Literature Review}

The presence of the continuous spectrum in $L^2(\Gamma \backslash G)$ always demands delicate analysis and one of the most promising tool to use is Selberg trace formula. The main obstruction is that for any $f \in C_c^{\infty}(G)$, the convolution operator $R_{\Gamma}(f)$ is not trace class if $\Gamma$ is not cocompact. In \cite{Mul} Muller proved that however the restriction of $R_{\Gamma}(f)$ on $L^2_d(\Gamma \backslash G)$ is trace class. The convergence of the trace of $R_{\Gamma}(f)|_{L^2_d(\Gamma \backslash G)}$ is due to a control on the growth of the continuous spectrum by certain intertwining operators on principal series representations (see Sec. \ref{sec:- simple trace formula}\ ).

To get rid of continuous spectrum one often tries to obtain an ample source of test functions $f$ such that the convolution operator $R_{\Gamma}(f)$ is supported on discrete spectrum. In \cite{LV} such a class of compactly supported test functions was constructed so that $R_{\Gamma}(f)$ is supported on the spherical irreducible representations (having nonzero $K$-fixed vectors); it has been refined recently in \cite{AM}.  In the latest breakthrough by \cite{BLZZ}, for a prescribed irreducible representation $\pi$ occurring in the discrete spectrum of $L^2(\Gamma \backslash G)$ with congruence lattice $\Gamma$ a source of test functions $f$ (but merely Schwartz class functions) can be found so that the convolution operator $R_{\Gamma}(f)$ is supported \textit{nearly} on the $\pi$-isotypic component inside the discrete spectrum. Consequently one can write the corresponding trace formula with quite simple spectral side. The cost against this simplicity of spectral side is that there is no flexibility on the support of the support of the test functions (as those are Schwartz class) which causes problem while comparing trace formulas for two lattices. 

\subsection{Methodology}
Let us now describe a few crucial points of the idea behind the proof of the Thm. \ref{main thm}. By condition (1), in the comparison of the trace formulas for $\Gamma_1$ and $\Gamma_2$ the continuous spectra cancel out. But on the geometric side the unipotent orbital integral still survives. 

\noindent Now the set of regular elements in $G$ (denoted by $G_{\text{reg}}$) is an open and dense subset of $G$. Therefore we can take our test functions to be supported in $G_{\text{reg}}$ which guarantees the vanishing of the unipotent orbital integral. Therefore we get a simple trace formula consisting of discrete spectral side and regular geometric side (like the cocompact case). Restricting the support of the test functions inside  $G_{\text{reg}}$ does not change the trace of $\pi(f)$ for any irreducible representation $\pi$ because of Thm. \ref{thmHC2} together with the fact that $G \setminus G_{\text{reg}}$ is of measure zero. The rest of the argument goes same as in \cite{BR}.

\subsection{Acknowledgement} The authors sincerely thank Prof. C. S. Rajan for many fruitful discussions. K. Mondal is supported partially by Prime Minister Research Fellowship (PMRF), Govt. of India.

\section{Preliminaries}

\subsection{Associate parabolics in $G$}\label{Associate parabolics in $G$} Let $G$ be a semisimple Lie group of rank one and $K$ be a 
maximal compact subgroup in $G$ which we fix henceforth. In this situation any proper parabolic subgroup $P$ of $G$ is a minimal parabolic. Moreover, any two parabolic subgroups are pairwise associate, i.e., $G$-conjugate. \medskip

Let $P = LN$ be the Levi decomposition of a given proper parabolic $P$ in $G$, where $L$ is the Levi part of $P$ and $N$ is the unipotent radical of $P$. We also have Langlands decomposition $P = MAN$ where $A$ is the split torus inside $L$ and $M: = Z_K(A)$ is the centraliser of $A$ in $K$.
 In particular, $M$ is a semisimple (compact) subgroup of $K$. Let $\Sigma=\Sigma(G,A)$ be the set of roots for the pair $(G, A)$. If $\Sigma^+$ is a positive system relative to the unipotent radical $N$ then $\Sigma^+$ consists of at most two roots $\{ \alpha , 2 \alpha\}$. The rootspaces of $\alpha$ and $2 \alpha$ are denoted by $\mathfrak{n}_{\alpha}$ and $\mathfrak{n}_{2\alpha}$ respectively and the corresponding Lie groups of are denoted by $N_{\alpha},N_{2 \alpha}$ respectively. So $\mathfrak{n}_{\alpha} \oplus \mathfrak{n}_{2\alpha}$ is the Lie algebra of $N$. Let $m_1=\dim \mathfrak{n}_{\alpha}$ and $m_2 = \dim \mathfrak{n}_{2 \alpha}$.  
\medskip

\subsection{$\Gamma$-cuspidal parabolics} 
Let $\Gamma$ be a lattice in $G$ such that $\text{vol}\ (\Gamma \backslash G) <  \infty$ (finite covolume). A proper parabolic subgroup $P=MAN$ is called a $\Gamma$-cuspidal parabolic subgroup of $G$ if $P \cap \Gamma \subset MN$; $N / N \cap \Gamma $ is compact and $MN / MN \cap \Gamma$ is compact. Let $E(G, \Gamma)$ be the set of all $\Gamma$-cuspidal parabolic subgroups in $G$. The set $E(G, \Gamma)$ breaks as a disjoint union of finitely many orbits under the conjugation action of $\Gamma$. Let us fix the representatives $P_1, P_2, \dots, P_{n_\Gamma}$ of these orbits. They parametrise the set of all cusps of $\Gamma$. Call $P_1=P$. Then there exists $k_i \in K$ with $k_1=1$ such that for each $i\in \{ 1,2, \dots, n\}$ and $P_i=M_iA_iN_i$, $M_i=k_i^{-1}Mk_i$, $A_i=k_i^{-1} A k_i$ and $N_i=k_i^{-1} N k_i$. Note that if $\Gamma$ is a co-compact lattice in $G$ then $E(G, \Gamma)$ is empty; equivalently there is no cusp.

\subsection{Conjugacy classes of hyperbolic/elliptic elements} 
In this section, we mention some results about the geometry of conjugacy classes of hyperbolic/elliptic elements in a lattice in a semisimple Lie group.\medskip

\begin{lemma}\label{lem1}
    For any hyperbolic element $\gamma \in \Gamma$ the $G$-conjugacy class $[\gamma]_G$ is a closed set of measure zero.
\end{lemma}

\begin{proof}
    Any hyperbolic element is semisimple and therefore its conjugacy class must be closed. To prove that $[\gamma]_G$ is a set of measure zero, it is enough to show that the conjugacy class $[\gamma]_G$ endows a submanifold structure in $G$ with lower dimension. Via the pullback of the projection map, $[\gamma]_G$ is homeomorphic to $G/G_{\gamma}$ and $G_{\gamma}$ contains a Cartan subgroup, where $G_{\gamma}$ is a centraliser of $\gamma$ in $G$.
\end{proof}

\begin{lemma}\label{lem2}
    Let $E$ be the union of all conjugacy classes $[\gamma]_G$ of all hyperbolic elements $\gamma \in \Gamma$. Then $E$ is a closed set of measure zero.
\end{lemma}

\begin{proof}
    Let $C$ be a compact subset of $G$. Since $E$ is a countable union of closed set $E \cap C$ is a finite union of closed sets. Therefore $E \cap C$ is closed. This holds for all compact subset $C$. So $E$ is closed.
\end{proof}

\subsection{Spectral Decomposition}\label{sec: - Spectral Decomposition}

For the parabolic subgroup $P=MAN$ let $\Gamma_M=\Gamma \cap M$. So, $\Gamma_M$ is a finite subgroup of the compact group $M$. Let $\widehat{M}$ be the set of irreducible representations of $M$. Consider the right regular representation of $M$ in the Hilbert space $L^2(\Gamma_M \backslash M)$. We have the Hilbert space decomposition $$L^2(M / \Gamma_M) \cong \widehat{\bigoplus_{\sigma \in \widehat{M}}} m(\sigma , \Gamma_M)\ \sigma.$$

Let $W=N_K(A) / M$ be the Weyl group. Since $G$ is of real rank one (i.e. $\text{dim}\ A =1$) we have $|W|=2$. Let $w$ be the nontrivial element of $W$. Let $\mathfrak{a}$ be the Lie algebra of $A$. An irreducible representation $\sigma \in \widehat{M}$ is called ramified if $w \sigma =\sigma$. Otherwise, $\sigma$ is called unramified.\smallskip

\noindent For any $\sigma \in \widehat{M}$, $\nu \in \mathfrak{a}^*_{\C}$, let $\pi_{\sigma, \nu}$ denote the principal series representation of $G$. We write, 

$$\pi(\sigma , \nu) =
\begin{cases}
  m(\sigma , \Gamma_M) \, \pi_{\sigma , \nu} & \text{if } \sigma \text{ is ramified}, \\[5pt]
  m(\sigma , \Gamma_M) \, \pi_{\sigma, \nu} \oplus m(w\sigma , \Gamma_M) \, \pi_{w\sigma , \nu} & \text{if } \sigma \text{ is unramified}.
\end{cases}$$\smallskip

\noindent Let $H(\sigma , \nu)$ be the representation space of $\pi(\sigma , \nu)$. Let $\mathfrak{a}_i$ be the Lie algebra of $A_i$. Via conjugation by appropriate $k_i$, each $\sigma \in \widehat{M}$ gives rise to a representation $\sigma_i \in \widehat{M_i}$ and each $\nu \in \mathfrak{a}^*_{\C}$ gives rise to $\nu_i \in \mathfrak{a}^*_{i, \C}$. For each $\sigma \in \widehat{M}$ and $\nu \in \mathfrak{a}^*_{\C}$ we put 
$$\pi_{\Gamma}(\sigma, \nu) = \bigoplus_{i=1}^{n} \pi(\sigma_i , \nu_i).$$
\noindent We denote the corresponding representation space by $H_{\Gamma}(\sigma , \nu)$.\smallskip 

\noindent The right regular representation $L^2(\Gamma \backslash G)$ of $G$ is not semisimple if $\Gamma$ is not co-compact. Let $L^2_d(\Gamma \backslash G)$ be the maximal semisimple subrepresentation of $L^2(\Gamma \backslash G)$. The orthogonal complement of $L^2_d(\Gamma \backslash G)$ in $L^2(\Gamma \backslash G)$ is the continuous spectrum denoted by $L^2_c(\Gamma \backslash G)$. The Hilbert direct integral decomposition of $L^2_c(\Gamma \backslash G)$ can be written as 

$$ L^2_c(\Gamma \backslash G) \cong \bigoplus_{\sigma \in \widehat{M}} \int_{i \R^+}^{\oplus} \pi_{\Gamma}(\sigma, \nu)\ d\nu .$$

\noindent Note that a lattice $\Gamma$ in $G$ is not co-compact if and only if the continuous spectrum $L^2_c(\Gamma \backslash G) \neq 0$.

\subsection{Results on Harish-Chandra characters}

We will be using the following two important results about Harish-Chandra characters of irreducible unitary representations of semisimple Lie groups (see \cite{Kn}). \medskip

\noindent\textbf{Harish-Chandra character distribution:} Let $(\pi, W_{\pi})$ be an irreducible unitary representation of $G$. For $f \in C_c^{\infty}(G)$ the convolution operator $\pi(f)$ on $W_{\pi}$ is trace class and its trace is denoted by $\chi_{\pi}(f):= \Tr (\pi(f))$. 

\begin{thm}\label{thmHC1}
    Let $\{\pi_i\}$ be a finite collection of mutually inequivalent irreducible representations of $G$. Then their distribution characters $\{\chi_{\pi_i}\}$ are linearly independent distributions on $G$.
\end{thm}

\begin{thm}\label{thmHC2}
    Let $\pi$ be an irreducible unitary representation of $G$. Then the distribution character $\chi_{\pi}$ is given by a locally integrable function $\phi_{\pi}$ on $G$ i.e. for $f \in C_c^{\infty}(G)$, 
    $$\chi_{\pi}(f)=\int\limits_G f(g)\ \phi_{\pi}(g)\ dg.$$
    Moreover, the restriction of $\phi_{\pi}$ to the set $G_{\rm reg}$ of regular elements is a real analytic function and is invariant under conjugation.
\end{thm}

\section{Simple Trace Formula for Arithmetic Lattices}\label{sec:- simple trace formula}
We describe the simple trace formula for the representation $L^2(\Gamma \backslash G)$ by \cite{War}. Let $C_c^{\infty}(G)$ be the algebra of all compactly supported smooth functions on $G$. Because of the presence of continuous spectrum in the right regular representation $R_{\Gamma}$ on $L^2(\Gamma \backslash G)$, the convolution operator $R_{\Gamma}(f)$ for $f \in C_c^{\infty}(G)$ is not a trace class operator in general. However, the restriction of $R_{\Gamma}(f)$ on the discrete part $L^2_d(\Gamma \backslash G)$ is a trace class operator. To describe the trace we need to set up a few notations. For $\gamma \in \Gamma$, let $\Gamma_\gamma$ and $G_\gamma$ denote the centraliser of $\gamma$ in $\Gamma$ and $G$, respectively. For any $\sigma \in \widehat{M}$ and $\nu \in \boldsymbol{i} \mathfrak{a}^*$ there exists a unitary intertwining operator $$C_{\Gamma \sigma}(\nu) : H_{\Gamma}(\sigma, \nu) \longrightarrow H_{\Gamma}(\sigma , -\nu).$$ 
\noindent Moreover, for any $K$-finite compactly supported smooth function $f$ on $G$ one can use the functional equation of Eisenstein series \cite{HC}, \cite{Lan} to derive 
$$ C_{\Gamma \sigma}(\nu)\ \pi_{\Gamma}(\sigma , \nu)(f) = \pi_{\Gamma}(\sigma , -\nu)(f)\ C_{\Gamma \sigma}(\nu).$$

\noindent In literature $C_{\Gamma \sigma}(\nu)$ is known as scattering matrix for $\Gamma$ associated with holonomy class $\sigma \in \widehat{M}.$\smallskip

\noindent Now we are ready to display the required trace formula.
\begin{equation} \label{trace formula} \Tr R_{\Gamma}(f)|_{L^2_d(\Gamma \backslash G)} = I_{\Gamma}(f) + H_{\Gamma}(f) + S_{\Gamma}(f)+ T_{\Gamma}(f) + U_{\Gamma}(f);
\end{equation}
\noindent where, the terms in the sum are given, respectively, by \smallskip

\begin{align*}
\text{(identity)} \quad I_\Gamma(f) &= \mathrm{Vol}(\Gamma \backslash G)\, f(1), \\
\text{(regular)} \quad H_\Gamma(f) &= \sum_{\{\gamma\} \text{: hyperbolic or elliptic}} \mathrm{Vol}(\Gamma_\gamma \backslash G_\gamma)\int_{G_\gamma \backslash G} f(x^{-1}\gamma x)\, d(G_\gamma x), \\
\end{align*}

\begin{align*}
\text{(scattering)} \quad S_\Gamma(f) &= \frac{1}{4\pi} \sum_{\sigma \in \widehat{M}} \int_{Re(\nu)=0} \mathrm{tr}\left(C_{\Gamma\sigma}(\nu)^{-1} \frac{d}{d\nu} C_{\Gamma\sigma}(\nu) \, \pi_\Gamma(\sigma,\nu)(f)\right)\, |d\nu|, \\
\text{(threshold)} \quad T_\Gamma(f) &= -\frac{1}{4} \sum_{\sigma \in \widehat{M}} \mathrm{tr}\left(C_{\Gamma\sigma}(0)\, \pi_\Gamma(\sigma,0)(f)\right), \\
\text{(unipotent)} \quad U_\Gamma(f) &= \frac{1}{|\alpha_1|} \left(C_1(\Gamma) T_1(f) + C_2(\Gamma) T_2(f) + C'_1(\Gamma) T'_1(f)\right).
\end{align*}\smallskip

\noindent Here, the definition of $T_i(f), T_1'(f)$ can be recalled from \cite{OsWar} as follows:

\begin{align*}
T_1(f)&=\dfrac{1}{A(\mathfrak{n}_{\alpha})} \int\limits_K \int\limits_N f(k^{-1} n k)\ dn\ dk,\\
 T_2(f)&= \frac{|\alpha_1|}{2} \left\{
\int_{G/G_{n_0}} f(x n_0 x^{-1}) \, d(G/G_{n_0})(x)
+
\int_{G/G_{n_0^{-1}}} f(x n_0^{-1} x^{-1}) \, d(G/G_{n_0^{-1}})(x)
\right\},\\
 T_1'(f) &= \frac{(m_1 + 2m_2)}{A(\mathfrak{n}_{\alpha})} \int_{A(n_1)} \int_{N_1} \int_{N_2} 
f(k^{-1} n_1 n_2 k) \cdot \log |\log(n_1)| \, dk \, dn_1 \, dn_2,
\end{align*}

\noindent where $n_0$ is a representative of the nontrivial unipotent orbit in $\mathfrak{n}_{2 \alpha}$, 
$A(\mathfrak{n}_{\alpha})$ denotes the area of the unit sphere in $\mathfrak{n}_{\alpha}$, 
and $G_n$ represents the centraliser of an element  $n$ in $G$.
\smallskip

\begin{remark}
    If $\Gamma$ is torsion free, then there is no contribution of elliptic orbital integral in the regular term. The constants $C_i(\Gamma), C_1'(\Gamma)$ that appear in the unipotent contribution were explicitly computed by DeGeorge \cite{DG}.
\end{remark}

For the purpose of this paper, we put a certain condition on the $K$-finite function $f \in C_c^{\infty}(G)$ which make the trace formula much more simpler. An element $x \in G$ is called regular if the centralizer $C_G(x)$ of $x$ in $G$ is a torus. Consider the set $G_{\rm reg}$ to be the set of all regular elements in $G$. Then $G_{\rm reg}$ is an open and dense subset of $G$. Let $f$ be a $K$-finite function in $C_c^{\infty}(G_{\rm reg})$. Since any element $n$ of a unipotent radical $N$ is not regular, the identity term and unipotent term in the above trace formula vanish for $f \in C_c^{\infty}(G_{\rm reg})$. Consequently we have, 

$$ \Tr R_{\Gamma}(f)|_{L^2_d(\Gamma \backslash G)} = H_{\Gamma}(f) + S_{\Gamma}(f) +T_{\Gamma}(f).$$

\noindent Moreover if, $$ L^2_{d}(\Gamma \backslash G) \cong \bigoplus_{\pi \in \widehat{G}} m(\pi, \Gamma)\ \pi,$$

\noindent then we have the simple trace formula $$\sum\limits_{\pi \in \widehat{G}} m(\pi, \Gamma) \Tr(\pi (f))=H_{\Gamma}(f) + S_{\Gamma}(f) +T_{\Gamma}(f).$$

\begin{remark}
    The scattering and threshold terms are essentially coming from the spectral side (more specifically given by the continuous spectrum). If $G$ has square-integrable representations then we take its matrix coefficient as the test function. In this case the scattering and threshold terms vanish. But then the sum in the spectral side is supported only at a finite set of square-integrable irreducible representations (for example, discrete series representations). 
\end{remark}

\vspace{10pt}

\section{Proof of the Thm. \ref{main thm} 
}
In this section, we compare the representation spectra for two lattices $\Gamma_1$ and $\Gamma_2$ in $G$ and prove Thm. \ref{main thm}.

\subsection{Equality of Continuous Spectrum} \label{subsec: Equality of continuous Spectrum}

\begin{defn}
Two lattices $\Gamma_1$ and $\Gamma_2$ of finite covolume in $G$ are said to be isocuspidal if there exists a finite set $\{ P_1, P_2, \dots, P_n\}$ of cuspidal parabolic subgroups which serves as a common set of representatives of cusps for both lattices $\Gamma_1$ and $\Gamma_2$.
\end{defn}

\noindent Note that this is automatic if both the lattices are co-compact. Moreover, it also holds if the lattices $\Gamma_1$ and $\Gamma_2$ have only one cusp.


\noindent Let us denote $P:=P_1$. Then $P$ admits a decomposition say, $MAN$. Let us denote $\Gamma_{i,M}=\Gamma_{i}\cap M$ for $i=1,2$. Then $\Gamma_{i,M}$ are finite subgroups of $M$. Consider the representation spaces $L^2(\Gamma_{i,M} \backslash M)$ of $M$. We further assume that they are isomorphic as representations of $M$ i.e., $L^2(\Gamma_{1,M} \backslash M) \cong L^2(\Gamma_{2,M} \backslash M)$. Note that if the lattices $\Gamma_i$ are torsion-free then this assumption holds. Furthermore, if $$L^2(\Gamma_{i,M} \backslash M) \cong \bigoplus\limits_{\sigma \in \widehat{M}} m(\sigma, \Gamma_{i,M}) \sigma \hspace{20pt} \text{for}\ i=1,2$$ and $m(\sigma, \Gamma_{1, M})=m(\sigma, \Gamma_{2,M})$ for all but finitely many $\sigma \in \widehat{M}$ then by \cite{LM} $m(\sigma, \Gamma_{1, M})=m(\sigma, \Gamma_{2,M})$ for all $\sigma \in \widehat{M}$ i.e. $L^2(\Gamma_{1,M} \backslash M) \cong L^2(\Gamma_{2,M} \backslash M)$.

By the definition of $H_{\Gamma}(\sigma , \nu)$ in Sec. 3, it is clear that if two lattices $\Gamma_1$ and $\Gamma_2$ satisfy the above two assumptions, then the representations $(\pi_{\Gamma_1}(\sigma , \nu), H_{\Gamma_1}(\sigma , \nu)))$ and $(\pi_{\Gamma_2} (\sigma , \nu), H_{\Gamma_2}(\sigma, \nu)))$ are isomorphic. Moreover, the intertwining operators $C_{\Gamma_1 \, \sigma}$ and $C_{\Gamma_2\, \sigma}$ are equal. Therefore, the scattering terms for the two lattices $\Gamma_1$ and $\Gamma_2$ in the trace formula are equal i.e. $S_{\Gamma_1}(f)=S_{\Gamma_2}(f)$ for all $f \in C_c^{\infty}(G)$. Similarly, $T_{\Gamma_1}(f)=T_{\Gamma_2}(f)$ for all $f \in C_c^{\infty}(G)$. \medskip

\subsection{Equality of Discrete Spectrum}\label{subsec:-Equality of Discrete Spectrum}

Now by the assumption on the lattices $\Gamma_1$ and $\Gamma_2$ and Subsec. 5.1, comparison of the trace formulas for the lattices $\Gamma_1$ and $\Gamma_2$ gives 
\begin{align*}
    \quad \Tr R_{\Gamma_1}(f)|_{L^2_d(\Gamma_1 \backslash G)} - \Tr R_{\Gamma_2}(f)|_{L^2_d(\Gamma_2 \backslash G)} &= H_{\Gamma_1}(f)-H_{\Gamma_2}(f) \hspace{20pt}\ \text{for all}\ f \in C_c^{\infty}(G_{\rm reg}). \\ 
\end{align*}

\noindent Let us denote the set of all $G$-conjugacy classes of hyperbolic and elliptic elements of $\Gamma_1$ and $\Gamma_2$ by $[\Gamma_1]_G$ and $[\Gamma_2]_G$, respectively. For a $G$-conjugacy class $[\gamma]_G$ of a hyperbolic or elliptic element $\gamma \in \Gamma_1 \cup \Gamma_2$ we denote $\text{vol}\ (\Gamma_{\gamma} \backslash G_{\gamma})$ by $a(\gamma, \Gamma)$. So we have, 

$$ \sum\limits_{\pi \in \mathcal{S}}t_{\pi} \chi_{\pi}(f)=\sum\limits_{[\gamma] \in [\Gamma_1]_G \cup [\Gamma_2]_G} (a(\gamma, \Gamma_1)-a(\gamma, \Gamma_2)) \int\limits_{G_{\gamma}\backslash G} f(x^{-1} \gamma x)\ d(G_{\gamma}x),$$

\noindent where, $t_{\pi}=m(\pi,\Gamma_1)-m(\pi, \Gamma_2)$, $\mathcal{S}$ is a finite set such that $t_\pi = 0$ for all $\pi \in \widehat G \setminus \mathcal{S}$ and $\chi_{\pi}(f)=\Tr (\pi(f))$ for all $f \in C_c^{\infty}(G_{\rm reg})$.\smallskip

\noindent Since $\mathcal{S}$ is a finite set, by Thm. \ref{thmHC2} there exists $\phi \in L^1_{loc}(G)$ such that $$\int\limits_G f(g)\ \phi(g)\ dg=\sum\limits_{[\gamma] \in [\Gamma_1]_G \cup [\Gamma_2]_G} (a(\gamma, \Gamma_1)-a(\gamma, \Gamma_2)) \int\limits_{G_{\gamma}\backslash G} f(x^{-1} \gamma x)\ d(G_{\gamma}x).$$\smallskip

\noindent Let $E$ be the union of all conjugacy classes of hyperbolic elements in $\Gamma_1 \cup \Gamma_2$. Then by Lem. \ref{lem2}, $E$ is a closed subset of $G$ with measure zero. Let $g \notin E$. There exists an open set $U \subseteq G_{\text{reg}}$ in  containing $g$ such that $U \cap E$ is empty. Therefore, for any $f \in C_c(U)$, we have 
$$\int\limits_G f(g)\ \phi(g)\ dg=0.$$

This means that $\phi$ is essentially $0$ on $U$. Note that $U$ is a neighbourhood of an arbitrary point $g$ away from $E$ and $E$ is a closed set of measure zero. Therefore, $\phi$ is essentially $0$ on $G_{\text{reg}}$. Finally, $\phi$ is essentially $0$ on $G$ because the set $G \setminus G_{\text{reg}}$ of non-regular points is a measure zero set. Hence, $$ \sum\limits_{\pi \in \mathcal{S}}t_{\pi} \chi_{\pi}(f)=0, \hspace{20pt} \text{for all}\ f \in C_c^{\infty}(G).$$

\noindent From the linear independence of characters (Thm. \ref{thmHC1}) we have $t_{\pi}=0$ for all $\pi \in \mathcal{S}$. Hence, $L^2_d(\Gamma_1 \backslash G) \cong L^2_d(\Gamma_2 \backslash G).$

\end{document}